\documentclass[12pt]{article}
\usepackage{amssymb}
\usepackage{amsmath}
\usepackage{amsbsy}
\usepackage{amsthm}
\usepackage{rotating}

\newcommand{\be}{\begin{equation}}
\newcommand{\ee}{\end{equation}}
\newcommand{\ba}{\begin{eqnarray}}
\newcommand{\ea}{\end{eqnarray}}
\newcommand{\ban}{\begin{eqnarray*}}
\newcommand{\ean}{\end{eqnarray*}}

\newtheorem{theo}{Theorem}[section]
\newtheorem{prop}[theo]{Proposition}
\newtheorem{lem}[theo]{Lemma}
\newtheorem{cor}[theo]{Corollary}
\newtheorem{defi}[theo]{Definition}
\newtheorem{con}[theo]{Conjecture}

\begin{document}

\title{Finite part of operator K-theory for groups finitely embeddable into Hilbert space and the  degree of non-rigidity of manifolds}
\author{Shmuel Weinberger and Guoliang Yu\footnote{The authors are partially supported by NSF.}}
\date{ }

\maketitle

{\noindent  Abstract:} In this paper, we study lower bounds on  the K-theory of the maximal $C^*$-algebra of  a discrete group based on the amount of torsion it contains. We call this the finite part of the operator K-theory and give a lower bound that is valid for a large class of groups, called the ``finitely embeddable groups''. The class of  finitely embeddable groups  includes all residually finite groups, amenable groups, Gromov's monster groups, virtually torsion free groups (e.g. $Out(F_n)$), and any group of analytic diffeomorphisms of an analytic connected manifold fixing a given point. It is  an open question if  every  countable group is finitely embeddable.
We apply this result to measure the degree of non-rigidity for any compact oriented manifold $M$ with dimension $4k-1$ $(k>1)$. We derive a lower bound on the rank of the structure group $S(M)$ in this case. For  a compact Riemannian manifold $M$ with dimension greater than or equal to $5$ and  positive scalar curvature metric, there is an abelian group $P(M)$ that measures the size of the space of all positive scalar curvature metrics on $M$. We obtain a lower bound on the rank of the abelian group
 $P(M)$ when  the compact smooth spin manifold $M$  has dimension $2k-1$ $(k>2)$  and the fundamental group of $M$ is finitely embeddable.

\section{Introduction}

The purpose of this paper is to find simple lower bounds for operator K-theory of discrete groups of the conceptual simplicity of the dimension of a vector space and use the lower bounds to study the degree of non-rigidity for compact manifolds and the size of the space of positive scalar curvature metrics for compact manifolds.
Our measure of simplicity will be that ``extending the ring'' should not lose information. In other words, we wish to have lower bounds for operator K-theory of discrete groups, that survive any inclusion into larger groups.

The classical Novikov conjecture uses group homology to give lower bounds, but the existence of acyclic groups kills such lower bounds. The Baum-Connes conjecture and Farrell-Jones conjecture give conjectured pictures of how K-theory and L-theory of groups  should be built out of  group homology and
 K-theory and L-theory of finite groups. The K-theory of finite groups is typically very complicated, but it seems that very little  survives inclusions--for example, the Whitehead group of symmetric groups is trivial \cite{O}.

However, the classical trace map on the group algebra can be used to detect $K_0$ elements of the group algebra coming from finite subgroups--for each finite order element $g$ in the group, the map summing up all coefficients of group elements in the conjugacy class of the element $g$  is a trace on the group algebra \cite{S}.
This shows that the number of conjugacy classes of the finite order elements in the group is a lower bound for $K_0$ of the real group algebra.
Our requirement that  the $K_0$ elements survive inclusion into larger groups means that the best we can reasonably hope for would be to have a copy of
${\mathbb Z}$ for every finite order that occurs among the elements of the group--since any two elements of the same order in a group are conjugate in an HNN extension containing that group.

In this paper, we describe a general class of groups for which we can prove such a lower bound for operator K-theory. This class of groups  include all residually finite groups, amenable groups, hyperbolic groups,
Burnside groups,
Gromov's monster groups, virtually torsion free groups (e.g.
$Out(F_n)$), and  any group of analytic diffeomorphisms of an analytic connected manifold fixing a given point.
We apply  our lower bound to measure the degree of topological  non-rigidity for compact oriented manifolds within a given homotopy type
 and the size of all positive scalar curvature metrics on compact spin manifolds.

Now, we describe our results in more detail.

Let $G$ be a countable group. An element $g\in G$ is said to have order $d$ if $ d$ is the smallest positive integer such that $g^d=e$, where $e$ is the identity element of $G$. If no such positive integer exists, we say that the order of $g$ is $\infty$.

If $g\in G$ is an element in $G$ with finite order $d$, then we can define an idempotent in the group algebra ${\mathbb  Q} G$ by:
$$p_g= \frac{1}{d}(\sum_{k=1}^{d} g^k).$$

For the rest of this paper, we denote the maximal group $C^*$-algebra of $G$  by $C^\ast(G)$. We define $ K_0 ^{fin}( C^*(G)) $,  the finite part of  $K_0 ( C^*(G))$, to be the abelian subgroup of $K_0 ( C^*(G))$ generated by $[p_g]$ for all elements $g\neq e$ in $G$ with finite order.

\begin{con} If $\{g_1, \cdots, g_n\}$ is a collection of elements in $G$ with distinct finite orders  such that $g_i\neq e$ for all $1\leq i\leq n$,
 then

\noindent (1) $\{[p_{g_1}], \cdots, [p_{g_n}]\}$ generates an abelian subgroup of $ K_0 ^{fin}( C^*(G)) $ with rank $n$;

\noindent (2) any nonzero element in the abelian subgroup of $ K_0 ^{fin}( C^*(G)) $  generated by $\{[p_{g_1}], \cdots, [p_{g_n}]\}$  is not in the image of the assembly map
 $\mu: K_0 ^G (EG) \rightarrow K_0 (C^*(G))$, where $EG$ is the universal space for proper and free $G$-action.
\end{con}

In fact, we can state a stronger conjecture in terms of K-theory elements coming from finite subgroups and the number of conjugacy classes of nontrivial finite order elements. Such a stronger conjecture follows from the strong Novikov conjecture but would not survive inclusion into large groups.

The following concept is due to Gromov.

\begin{defi}
A countable discrete group $G$ is said to be coarsely embeddable into Hilbert space $H$ if there exists a map
$f: G \rightarrow H$ satisfying

\noindent (1) for any finite subset $F\subseteq G$, there exists $R>0$ such that
if $g^{-1}h \in F$, then $d(f(g), f(h))\leq R$;

\noindent (2) for any $S>0$, there exists a finite subset $E\subseteq G$ such that if $g^{-1}h\in G-E$, then $d(f(g), f(h)) \geq S$.
\end{defi}

The class of groups coarsely embeddable into Hilbert space includes amenable groups \cite{BCV}, hyperbolic groups \cite{Sela}, and  linear groups \cite{GHW}.
However, Gromov's monster groups are not coarsely embeddable into Hilbert space \cite{G, AD}.
The importance of the concept of coarse embeddability is due to the theorem that the strong Novikov conjecture holds for groups coarsely embeddable into Hilbert space \cite{Y1, STY}. In \cite{KY}, Kasparov and Yu introduced a weaker condition, coarse embeddability into Banach spaces with Property H, and
proved the strong Novikov conjecture for groups coarsely embeddable into Banach spaces with Property H.

The following concept is   more flexible  than coarse embeddability into Hilbert space.

\begin{defi}  A countable discrete group $G$ is said to be finitely embeddable into Hilbert space $H$ if for any finite subset $F\subseteq G$,
there exists a group $G'$ coarsely embeddable into $H$ such that there is a map $\phi: F\rightarrow G'$ satisfying

\noindent (1) $\phi(gh)=\phi(g) \phi (h)$ if $g, h\in F$ and $gh\in F$;

\noindent (2) if $g$ is a finite order  element in $F$, then $order(\phi(g))= order (g)$.

\end{defi}

We mention that the class of groups finitely embeddable into Hilbert space include all residually finite groups, amenable groups, hyperbolic groups,
Burnside groups,
Gromov's monster groups, virtually torsion free groups (e.g.
$Out(F_n)$), and  any group of analytic diffeomorphisms of an analytic connected manifold fixing a given point. It is an open question to construct a group  not finitely embeddable into Hilbert space. We can similarly define a concept of finite embeddability into Banach spaces with Property H.

The general validity of Conjecture 1.1 is still open.
The following result proves Conjecture 1.1 for a large class of groups.

\begin{theo} Conjecture 1.1 holds for groups finitely embeddable into Hilbert space.
\end{theo}

We define $N_{fin} (G)$ to be the cardinality of the following subset of positive integers:
$$\{ d: \exists ~~g\in G~~ s. ~t.~~ g\neq e, ~~order(g)=d\}.$$

If $M$ is a compact oriented manifold,  the structure group $S(M)$ is the abelian group of equivalence classes of all pairs  $(f, M')$ such that
$M'$ is a  compact oriented  manifold and
$f: M' \rightarrow M,$ is an orientation preserving homotopy equivalence \cite{R}.
The rank of $S(M)$ measures the degree of non-rigidity for $M$.

The following result explains why it is interesting  to study the finite part of  $K_0 ( C^*(G))$.

\begin{theo}
Let $M$ be a compact oriented manifold with dimension $4k-1$  $(k>1)$ and $\pi_1(M)=G$.
If Conjecture 1.1 holds for $G$, then the rank of the structure group $S(M)$ is greater than or equal to $N_{fin}(G)$.
\end{theo}

The following result is a consequence of Theorems 1.4 and 1.5.

\begin{cor} Let $M$ be a compact oriented manifold with dimension $4k-1$  $(k>1)$ and $\pi_1(M)=G$.
If $G$ is finitely embeddable into Hilbert space, then the rank of the structure group $S(M)$ is greater than or equal to $N_{fin}(G)$.
\end{cor}

We will also prove that, for a restricted class of fundamental groups, the different elements in the structure  set detected by our method represent different manifolds.

Let  $M$ be a compact smooth manifold with dimension greater than or equal to $5$. If $M$ has a metric with positive scalar curvature, following Stolz \cite{S2, RS},
we will introduce an abelian group $P(M)$ of   concordance classes of  all positive scalar curvature metrics on $M$ (cf. Section 4 for a more precise definition).
The following theorem gives an estimation on the ``size'' of the space of positive scalar curvature metrics on $M$ when the fundamental group is finitely embeddable into Hilbert space.
\begin{theo}
\noindent (1) Let  $M$ be a compact smooth spin manifold with dimension $2k-1$ $(k>2)$. If $M$ has a metric with positive scalar curvature and $\pi_1(M)=G$ is finitely embeddable into Hilbert space, then the rank of the abelian group $P(M)$ is greater than or equal to $N_{fin}(G)$;

\noindent (2) Let  $M$ be a compact smooth spin manifold with dimension $4k-1$ $(k>1)$. If $M$ has a metric with positive scalar curvature and $\pi_1(M)=G$ is finitely embeddable into Hilbert space, then the rank of the abelian group $P(M)$ is greater than or equal to $N_{fin}(G)$+1.
\end{theo}

We remark that the main results in this paper remain to be true under the weaker condition of finite embeddability into Banach spaces with Property H.

In general, we also have the following results on the ``size'' of the  space of positive scalar curvature metrics on $M$.

\begin{theo}
\noindent (1) Let $M$ be  a compact smooth spin manifold with a positive scalar curvature metric and  dimension  $2k-1$ ($k>2$).
 If $\pi_1(M)$ is not torsion free, then the rank of the abelian group $P(M)$ is greater than or equal to one;

\noindent (2) Let $M$ be  a compact smooth spin manifold with a positive scalar curvature metric and dimension  $4k-1$ ($k>1$).
 If $\pi_1(M)$ is not torsion free, then the rank of the abelian group $P(M)$ is greater than or equal to two.
\end{theo}

Part (1) of the above theorem is motivated by  and is a slight refinement of a theorem of Piazza-Schick stating that if $M$ is a compact smooth spin manifold with a positive scalar curvature metric and  dimension  $2k-1$ ($k>2$), and $\pi_1(M)$ is not torsion free,
than the moduli space of positive scalar curvature metrics on $M$ has infinitely many connected components \cite{PS}. In \cite{GL}, Gromov-Lawson proved that the space
of positive scalar curvature metrics on a compact spin manifold $M$ has infinitely many connected components if the dimension of $M$ is $4k$ ($k>1$).

This paper is organized as follows. In Section 2, we study the finite part of K-theory for group $C^*$-algebras.
In Section 3, we discuss its application to non-rigidity of compact oriented  manifolds. In Section 4, we discuss its application to positive scalar metrics.
In Section 5, we study groups finitely embeddable into Hilbert space.

The authors wish to thank Sherry Gong, Jonathan Rosenberg, Melissa Liu, Rufus Willett, Zhizhang Xie, and Romain Tessera  for very helpful discussions. This paper is written during the second author's stay at the Shanghai Center for Mathematical Sciences (SCMS). The second author wishes to thank SCMS for providing a wonderful environment.

\section{Finite part of K-theory for group $C^*$-algebras}

In this section, we study the finite part of K-theory for group $C^\ast$-algebras. In particular,
 we introduce a concept of finite representability and prove  that Conjecture 1.1 is closed under finite representability.

\begin{defi} Let $\mathcal{F}$ be a family of countable groups. We say that  a countable group $G$ is finitely representable in
$\mathcal{F}$ if for any finite subset $F\subseteq G$,
there exists a group $G'\in  \mathcal{F}$ such that there is a map $\phi: F\rightarrow G'$ satisfying

\noindent (1) $\phi(gh)=\phi(g) \phi (h)$ if $g, h\in F$ and $gh\in F$;

\noindent (2) if $g$ is a finite order  element in $F$, then $order(\phi(g))= order(g)$.
\end{defi}

\begin{theo} Let $\mathcal{F}$ be a family of countable groups for which Conjecture 1.1 holds.
If $G$ is finitely representable in $\mathcal{F}$, then Conjecture 1.1 holds for $G$.
\end{theo}
\proof
Let $\{F_k\}_{k=1}^\infty$ be an increasing family of finite subsets of $G$ such that

\noindent (1) $\cup_{k=1}^\infty F_k =G;$

\noindent (2) for each $k$, there exists a group $G_k'$ in $\mathcal{F}$ and a map $f_k: F_k \rightarrow G_k'$ satisfying
$f_k(gh)=f_k (g) f_k (h)$ for all $g, h\in F_k$ such that $gh\in F_k$;

\noindent (3) if $g$ is a finite order element   in $F_k$, then $order(f_k(g)) =order (g).$

Let $\phi_k$ be the linear map: $ \mathbb{C}F_k \rightarrow C^* (G_k')$, extending $f_k$.
The family of maps $\{\phi_k\}_{k=1}^\infty$ induces a $\ast$-homomorphism:
$$\phi:  \mathbb{C} G \rightarrow \prod_{k=1}^\infty C^* (G_k')/(\oplus_{k=1}^{\infty} C^* (G_k'))$$
defined by $$\phi(a)=[\prod_{k=n_a}^{\infty} \phi_k(a)]$$ for all $a\in  \mathbb{C} G ,$ where $n_a$ is a positive integer such that
the support of $a$ is contained in $F_k$ for all $k\geq n_a$.

By the definition of the maximal group $C^*$-algebra, $\phi$ can be extended to a $\ast$-homomorphism
(still denoted by $\phi$):
$$\phi:  C^* G \rightarrow \prod_{k=1}^\infty C^* (G_k')/(\oplus_{k=1}^{\infty} C^* (G_k')).$$

We have the following six term exact sequence:

{\tiny
$$\begin{array}[c]{ccccc}
K_0 (\oplus_{k=1}^{\infty} C^* (G_k')) &\stackrel{i_*}{\rightarrow}& K_0 (\prod_{k=1}^\infty C^* (G_k'))&\stackrel{\pi_*}{\rightarrow}& K_0( \prod_{k=1}^\infty C^* (G_k')/(\oplus_{k=1}^{\infty} C^* (G_k')     )) \\
\uparrow& & & &\downarrow\scriptstyle{\partial}\\
K_1 ( \prod_{k=1}^\infty C^* (G_k')/(\oplus_{k=1}^{\infty} C^* (G_k') )  )&\leftarrow& K_1 ( \prod_{k=1}^\infty C^* (G_k'))&\stackrel{j_*}{\leftarrow} & K_1(\oplus_{k=1}^{\infty} C^* (G_k') )
\end{array}.$$ }

Observe that $j_*$ is injective. By exactness, $\pi_*$ is surjective.
If $g$ is an element in $G$ with finite order $d$, there exists a positive integer  $N_g>0$ such that $g^m\in F_{N_g} $ for all $m$.
Let $p_k =\phi_k (p_g)$ for all $k\geq N_g$. Define $$p'_g=\prod_{k=N_g}^{\infty} p_k.$$

Let  $\{g_1, \cdots, g_n\}$ be a collection of elements in $G$ with $order(g_l)=d_l$ such that $g_l\neq e$ for each $l$ and $order(g_{l_1})< order(g_{l_2})$ when $1\leq l_1<l_2\leq n$.
The rank of the subgroup of $K_0(C^\ast(G))$ generated by $\{[p_{g_1}], \cdots, [p_{g_n}]\}$  is greater than or equal to
 the rank of the abelian subgroup of \linebreak $K_0( \prod_{k=1}^\infty C^* (G_k')/(\oplus_{k=1}^{\infty} C^* (G_k')     ))$
 generated by $\{[ \phi(p_{g_1})],  \cdots, [ \phi(p_{g_n})]\} $.
 In the same time, the rank of  the abelian subgroup of $K_0( \prod_{k=1}^\infty C^* (G_k')/(\oplus_{k=1}^{\infty} C^* (G_k')     ))$
 generated by $\{[ \phi(p_{g_1})],  \cdots, [ \phi(p_{g_n})]\} $   is the same as the rank
 of the abelian subgroup of $K_0 (\prod_{k=1}^{\infty} C^* (G_k'))$ generated by $\{[p'_{g_1}], \cdots, [p'_{g_n}]\}$. This is because, if $z$ is a nonzero element in the abelian subgroup of $K_0 ( \prod_{k=1}^{\infty} C^* (G_k'))$ generated by $\{[p'_{g_1}], \cdots, [p'_{g_n}]\}$, then
 $z$ is not in the image of $i_*$.
 As a consequence, we have
  $$\pi_* (z) \neq 0 ~~in~~ K_0( \prod_{k=1}^\infty C^* (G_k')/(\oplus_{k=1}^{\infty} C^* (G_k')     )) .$$
 By assumption, the rank
 of the abelian subgroup of $K_0 (\prod_{k=1}^{\infty} C^* (G_k'))$ generated by $\{[p'_{g_1}], \cdots, [p'_{g_n}]\}$ is $n$.
 Hence the rank of  the abelian subgroup of $K_0( \prod_{k=1}^\infty C^* (G_k')/(\oplus_{k=1}^{\infty} C^* (G_k')     ))$
 generated by $\{[ \phi(p_{g_1})],  \cdots, [ \phi(p_{g_n})]\} $   is $n$. It follows that the rank of the subgroup of $K_0(C^\ast(G))$ generated by $\{[p_{g_1}], \cdots, [p_{g_n}]\}$  is also $n$.

 Let $x$ be a nonzero element in the abelian subgroup of $K_0^{fin} (C^*(G))$ generated by $\{[p_{g_1}], \cdots, [p_{g_n}]\}$.
 Assume by contradiction that $x$ is in the image of the map $\mu: K_0^G(EG) \rightarrow K_0(C^*(G)).$
 Let $N$ be a positive integer such that $g_i ^k \in F_N$ for all $i$ and $k$.
 This implies that the element $\prod_{k=N}^\infty \phi_\ast (x)$ in $ \prod_{k=1}^\infty K_0 (C^* (G_k'))/(\oplus_{k=1}^{\infty} K_0(C^* (G_k')     )) $
 is in the image of the map $[\prod_{k=1}^\infty \mu]:$
 $$ \prod_{k=1}^\infty K_0^{G_k'} (EG_k')/(\oplus_{k=1}^{\infty} K_0^{G_k'}(EG_k'     ))  \rightarrow \prod_{k=1}^\infty K_0 (C^* (G_k'))/(\oplus_{k=1}^{\infty} K_0(C^* (G_k')     ))
 .$$ This implies that $\phi_\ast (x)$ is in the image of $\mu: K_0^{G_k'} (EG_k')\rightarrow K_0(C^* (G_k')$ for some large $k\geq N$.
 This is a contradiction with the assumption that $G_k'$ satisfies Conjecture 1.1.
 \qed

We remark that the above result is still open for reduced group $C^\ast$-algebras.

\begin{theo} (1) If there exists an element $g \neq e$ in $G$ with finite order, then the rank of the abelian subgroup of  $K_0 (C^\ast(G))$ generated by $[p_g]$ is one and any nonzero element in the abelian subgroup is not in the image of the map $\mu:
K_0^G(EG)\rightarrow K_0(C^*(G))$;
(2) If there exist two elements  in $G$ with finite orders
such that the order of $g_1$ is not equal to the order of $g_2$, then the rank of abelian subgroup of $K_0 (C^\ast(G))$ generated by $[p_{g_1}]$
and $[p_{g_2}]$ is two.
\end{theo}
\proof
(1) Let $tr: \mathbb{C}G\rightarrow \mathbb{C}$ be the canonical trace defined by:
$$tr(\sum_{g\in G} c_g g)=c_e.$$
It is easy to see that $tr$ extends to a finite trace on both $C^*(G)$ and $C^*_r (G)$.
Let $d$ be the order of $g$.
We have $$tr(p_g)=\frac{1}{d}.$$ This shows that the abelian subgroup generated by $[p_g]$ in $K_0 (C^\ast(G))$ or
$K_0( C^\ast_r(G))$ has rank one.

Let $\tau$ be the trace on $\mathbb{C}G$ defined by: $$\tau(\sum_{g\in G} c_g g)= \sum_{g\in G} c_g.$$
$\tau$ is a $\ast$-homomorphism from  $\mathbb{C}G$ to $\mathbb{C}.$ By the definition of $C^\ast(G)$, $\tau$ extends to a $\ast$-homomorphism
from $C^*(G)$ to $\mathbb{C}.$ It follows that $\tau$ is a trace on $C^*(G)$.

Assume by contradiction that a nonzero element $x$ in the abelian subgroup is  in the image of the map $\mu:
K_0^G(EG)\rightarrow K_0(C^*(G))$. This means that there exists $y\in K_0^G(EG)$ satisfying $\mu(y)=x$. By Atiyah's $L^2$-index theorem in \cite{A} and
 Theorem 6.1 in the Appendix of this paper, we have $$tr(x) =\tau(x)=index (y),$$
where $index(y)$ is the Fredholm index of K-homology class $y$. This is a contradiction with
the fact that $tr (p_g)-\tau(p_g)\neq 0$ and $x$ is a nonzero multiple of $[p_g]$.

(2)
We have $$tr(p_{g_1})=\frac{1}{d_1}, ~~~~tr(p_{g_2})=\frac{1}{d_2};$$
$$\tau(p_{g_1})=\tau(p_{g_2})=1.$$
This proves that the abelian subgroup generated by $[p_{g_1}]$ and $[p_{g_2}]$  in $K_0 (C^\ast(G))$ has rank two.
\qed

We remark that the first part of  (1) in the above theorem holds for both maximal and reduced group $C^\ast$-algebras.  It remains open to prove the reduced group
$C^\ast$-algebra analogue of second part of (1) in the above theorem and the reduced group
$C^\ast$-algebra analogue of part (2) in the above theorem. It is also an  open question whether any nonzero element in the abelian subgroup of $K_0 (C^\ast(G))$ generated by  $[p_{g_1}]$
and $[p_{g_2}]$ is not in the image of the map $\mu:
K_0^G(EG)\rightarrow K_0(C^*(G))$ for any two nontrivial finite order elements $g_1$ and $g_2$ in $G$.

\section{Applications to degree of non-rigidity of manifolds}

In this section, we use   the finite part of K-theory for maximal group $C^*$-algebras to estimate the degree of
non-rigidity. We also estimate the rank of the finite part of K-theory for maximal group $C^*$-algebras in terms of the size of the set of finite order elements
when the groups are finitely embeddable into Hilbert space.

Let $G$ be a countable group. Let $X$ be a locally compact space with a proper and free cocompact action of $G$.

Let $C_0(X)$ be the algebra of all continuous functions on $X$ vanishing at infinity.
If $H$ is a $X$-module (i.e. $H$ is a Hilbert space with an action of $C_0(X)$), we define the support of an operator $T: H\rightarrow H$, $support(T)$, to be the complement of the
subset of $X\times X$ consisting of all points $(x,y)\in X\times X$ such that there exists $f$ and $g$ in $C_0(X)$ such  that
$f(x)\neq 0$, $g(y)\neq 0$, and $fTg=0.$ In this paper, all $X$-modules are ample in the sense that no non-zero function of $C_0(X)$ acts
the $X$-modules as a compact operator.

We  endow $X$ with a $G$-invariant proper metric
$d$ (compatible with the topology of $X$).  An operator on $H$ is said to have finite propagation \cite{Roe}
if $$ propagation(T)=sup \{ d(x,y): (x, y)\in ~~support(T)\}<\infty.$$

Recall that an operator  $T$ on $H$ is called locally traceable if for any pair of compactly supported functions $f$ and $g$,
$fTg$ is a trace class operator. For each $p\geq 1$, we can similarly define the concept of locally Schatten-p class operators.
Let $({\cal S}_1 X)^G$ be the algebra of $G$-invariant and locally traceable operators
on $H$ with finite propagation. We observe that $({\cal S}_1 X)^G$ is isomorphic to ${\cal S}_1 G$, the group algebra over the ring of $ {\cal S}_1$.
There is a canonical trace $tr$ on ${\cal S}_1 G$:
$$tr (\sum_{g\in G} s_g g)= trace (s_e).$$
We use the same notation $tr$ to be the corresponding trace on $({\cal S}_1 X)^G$.

\begin{lem}
Let $H$ be as above, let $T \in ({\cal S}_1 X)^G$.
 If there exists $\delta>0$ such that $$support (T)\subseteq \{(x,y)\in X\times X: d(x, y)\geq \delta\},$$
then $tr(T)=0.$
\end{lem}
\proof Let $Z\subseteq X$ be a Borel fundamental domain of $X$, i.e.
(1) $\cup_{g\in G} ~~gZ= X$; (2) $gZ\cap Z =\emptyset$ for any $g\neq e$;
(3) $Z$ is bounded.

Let $\chi_Z$ be the characteristic function of $Z$. We have $$tr (T) =trace(\chi_Z T\chi_Z).$$
Note that the support of $ \chi_Z T\chi_Z $ is contained in $\{ (x, y)\in X\times X: d(x, y)\geq \delta\}$.

We decompose $$Z=Z_1\cup \cdots \cup Z_N$$ such that $$diameter (Z_i) <\frac{\delta}{2}.$$
Let $H_Z=\chi_Z H, H_i=\chi_{Z_i} H.$  Let $T_{i,j} $ be the restriction of $\chi_{Z_i}T\chi_{Z_j}$ to $H_j$ and let $T_Z$ be the restriction of
$\chi_Z T\chi_Z$ to $H_Z$.

We have a Hilbert space decomposition: $$H_Z=\oplus_{i=1}^N H_i.$$ We have a corresponding matrix representation of $T_Z$:
$$T_Z= (T_{i,j}) _{1\leq i, j\leq N} .$$ By the support condition on $T_Z$, we have $T_{i,i}=0$ for all $i$.
This implies $ trace(T_Z)=0.$
\qed

Let $G$ be a countable group and let ${\cal S}$ be the ring of Schatten class operators.
Let ${\cal S}G$ be the group algebra over the ring ${\cal S}$ \cite{Y2}. Let $j: \mathbb{C}G\rightarrow {\cal S} G$ be inclusion homomorphism defined by:
$$j(a)= p_0 a $$ for all $a\in \mathbb{C}G$, where $p_0$ is a rank one projection in ${\cal S} .$
We define the finite part of $K_0 ({\cal S}G)$ to be the abelian subgroup of $K_0 ({\cal S}G)$ generated by $j_\ast [p_g]$ for all
finite order elements $g\neq e$ in $G$.

\begin{lem}
Any nonzero element in the finite part of $K_0 ({\cal S}G)$  is not in the image of the assembly map:
$$A: H_0^{Or~G}(EG, {\mathbb K}({\cal S})^{-\infty}) \rightarrow K_0 ({\cal S}G),$$ where the assembly map $A$ is defined as in \cite{BFJR}.
\end{lem}
\proof
Let $g$ be a finite order element $g\neq e$ in $G$.
Let $tr_g: \mathbb{C}G \rightarrow \mathbb{C}, $ be the trace defined by:
$$tr_g(\sum_{\gamma} c_{\gamma} \gamma) = \sum_{ \gamma \in C(g) } c_{\gamma},$$
where $C(g)$ is the conjugacy class of $g$, i.e. $$C(g)=\{ \gamma \in G:  \exists~h\in G: ~h^{-1}\gamma h=g \}.$$

For any integer $m>0$, let ${\cal S}_m$ be the ring of Schatten-$m$ class operators.
Let $n=2k$ be the smallest even number greater than or equal to $m$.
Define an $n$-cyclic cocycle $\tau_g^{(n)}$ on $ {\cal S}_m G$ by:
$$\tau_g^{(n)} ( a_0, a_1, \cdots, a_n) =  \sum_{ \gamma\in C(g) } tr(\gamma ^{-1}a_0a_1\cdots a_n)$$
for all $a_i \in {\cal S}_m G$, where $C(g)$ is the conjugacy class of $g$ and $tr: {\cal S}_1 G \rightarrow \mathbb{C},$
is the trace defined by: $$tr(\sum_{\gamma \in \Gamma} b_\gamma \gamma) =trace(b_e).$$

If  $S$ is the suspension operator in Connes' theory of cyclic cohomology theory \cite{C}, then   $$ (S^k tr_g)(a_0, \cdots, a_n) = tr_g(a_0 \cdots a_n)$$
for all $a_i \in {\mathbb C}G.$

We have $$S^k tr_g=j^\ast \tau_g^{(n)}$$ as cyclic cocycles over ${\mathbb C}G.$
 In particular, this implies that the pairing between $tr_g$ and any idempotent $p$ in the matrix algebra of ${\mathbb C} G$ is the same as the pairing between $\tau_g^{(n)}$ and $j(p)$, i.e. $$\tau_g^{(n)}(j(p), \cdots, j(p))= tr_g(p),$$ where
 the cyclic cocycle $\tau_g^{(n)}$ is naturally extended to a cyclic cocycle over the matrix algebra of ${\mathbb C} G.$

If an element $[z]$ is in the image of the assembly $$A: H_0^{Or~G}(EG, {\mathbb K}({\cal S})^{-\infty}) \rightarrow K_0 ({\cal S}G),$$
then there exist $m>0$ and  a locally compact simplicial complex $X$ with a proper and cocompact action of $G$ such that
$[z]$ is in the image of the assembly map $$A: H_0^{Or~G}(X, {\mathbb K}({\cal S}_m)^{-\infty}) \rightarrow K_0 (({\cal S}_m X)^G) \cong K_0 ({\cal S}_m G),$$
 where  $X$ is endowed  with a $G$-invariant proper metric
$d$ (compatible with the topology of $X$) and $({\cal S}_m X)^G$ is the algebra of $G$-invariant and locally Schatten-m  operators
on $H$ with finite propagation ($({\cal S}_m X)^G$ is isomorphic to $ {\cal S}_m G $).
By cocompactness of the $G$-action on $X$, there exists $\delta>0$ such that $d(\gamma x, x) \geq 10\delta$ for all $x\in X$ and $\gamma\neq e$ in $G$.
If $H$ is a $G$-$X$-module, let $({\cal S}_m X)^G$ be the algebra of $G$-invariant and locally Schatten-m class operators
on $H$ with finite propagation.
The class $[z]$ can be represented by $$[q]-  \bigl(\begin{smallmatrix}
1&0\\ 0&0
\end{smallmatrix} \bigr),$$ where $q$ is an idempotent in the matrix algebra of $(({\cal S}_m X)^G)^+$
(~$(({\cal S}_m X)^G)^+$ is obtained from $({\cal S}_m X)^G$ by adjoining a unit).
We naturally extend the  cyclic cocycle $\tau_g^{(n)}$ to a cyclic cocycle  over the matrix algebra of $(({\cal S}_m X)^G)^+$ (still denoted by $\tau_g^{(n)}$)
by setting  $\tau_g^{(n)}(a_0, \cdots, a_n)=0$
 if $a_i$ is a scalar matrix for some $i$.

By the definition of the assembly map in \cite{BFJR} and the fact that this assembly map coincides with the classic assembly map
 (Corollary 6.3 in \cite{BFJR}), we can write
$$q=q_1+ \bigl(\begin{smallmatrix}
1&0\\ 0&0
\end{smallmatrix} \bigr) $$ such that $q_1$ is an element in the matrix algebra over the algebra of $G$-invariant and locally Schatten-m class operators
on a $G$-$X$-module $H$ with finite propagation satisfying $$propagation(q_1)<\delta/(n+1).$$
We have $$propagation (q_1^{n+1})<\delta.$$ This, together with Lemma 3.1, implies  $$tr (\gamma^{-1} q_1^{n+1})=0$$
for all $\gamma \neq e$ in $G$.
It follows that $$ \tau_g^{(n)} (q, \cdots, q)=0.$$ Recall that the pairing between the cyclic cocycle $\tau^{(n)}_g$  and $K_0 (({\cal S}_m X)^G)$
gives a homomorphism $(\tau^{(n)}_g)_\ast: K_0 ( {\cal S}_m G ) \rightarrow {\mathbb C}$.
This implies that $$(\tau^{(n)}_g)_\ast (z)=0~~~~~~~~~ \cdots \cdots~~~~~~~~~ (\ast) $$ for any finite order element $g\neq e$ in $G$ and  any element $[z]$  in the image of the assembly map
$$A: H_0^{Or~G}(X, {\mathbb K}({\cal S}_m)^{-\infty}) \rightarrow K_0 (({\cal S}_m X)^G).$$

If  $x$ is a nonzero element in the finite part of $ K_0 ( {\cal S}G) $, then $x$ is a nonzero element in the finite part of
$ K_0 ( {\cal S}_mG) $ if  $m$ is some large natural number.  We can write $$z=\sum_{i=1}^l c_i [j(p_{g_i})]$$ for some $c_i \in {\mathbb Z}$
and finite order elements $g_i\neq e \in G$.
We can assume that the conjugacy class of $g_{i_1}$ is different from the conjugacy class of $g_{i_2}$ when $i_1\neq i_2$ (otherwise $[j(p_{g_{i_1}})]=[j(p_{g_{i_2}})]$ in
$ K_0 ( {\cal S}_m G) $). Without loss of generality, we can assume that $c_l\neq 0$  and $order(g_{i_1})\leq order(g_{i_2})$ when $i_1\leq i_2$. Since the conjugacy class of $g_{i_1}$ is
different from the conjugacy class of $g_{i_2}$ when $i_1\neq i_2$, we have
$$tr_{g_i} (p_{g_i})\neq 0,~~~~~~tr_{g_{i_2}} (p_{g_{i_1}})=0~~~~~when ~~~i_1<i_2.$$
The above facts, together with  the the identity
$$\tau^{(n)}_{g_l}(j(p_{g_i}), \cdots, j(p_{g_i})) = tr_{g_l} (p_{g_i}),$$
imply   that $$
(\tau^{(n)}_{g_l})_\ast (x) \neq 0.$$
Combining this with the  equation $(\ast)$, we obtain Lemma 3.2.
\qed

\begin{lem} If $\{g_1, \cdots, g_n\}$ is a collection of finite order elements in $G$  with distinct conjugacy classes, then $\{j_\ast[p_{g_1}], \cdots, j_\ast [p_{g_n}]\}$ generates an abelian subgroup of $ K_0 ( {\cal S}G) $ with rank $n$, where $j_\ast$ is as in Lemma 3.2.
In particular, if $\{g_1, \cdots, g_n\}$ is a collection of nontrivial finite order elements in $G$  with distinct orders, then $\{j_\ast[p_{g_1}], \cdots, j_\ast [p_{g_n}]\}$ generates an abelian subgroup of $ K_0 ( {\cal S}G) $ with rank $n$.
\end{lem}
\proof Lemma 3.3  essentially follows from the argument in the last paragraph of the proof of Lemma 3.2.

It suffices to prove that $\{j_\ast[p_{g_1}], \cdots, j_\ast [p_{g_n}]\}$ generates an abelian subgroup of $ K_0 ( {\cal S}_m G) $ with rank $n$
for any natural number $m\geq 1$.

Let $g$ be an element in $G$ with finite order. Let $tr_g$ be defined as in the proof of Lemma 3.2.
By  assumption, $\{g_1, \cdots, g_n\}$ has distinct conjugacy classes. Without loss of generality, we can assume that $order(g_{i_1})\leq order(g_{i_2})$ when $i_1\leq i_2$.
 As a consequence,
 we have
$$tr_{g_i}(p_{g_i})\neq 0, ~~~~~tr_{g_{i_2}} (p_{g_{i_1}})=0~~~~when ~~~i_1<i_2.$$

As in the proof of Lemma 3.2,  $tr_{g_i}$ can be lifted to a cyclic cocycle $\tau^{(n)}_{g_i}$ on $ {\cal S}_m G$ when $n=2k\geq m$.
Hence for each m, there exist homomorphisms $(\tau^{(n)}_{g_i})_\ast: K_0 ( {\cal S}_m G ) \rightarrow {\mathbb C}$ such that
$$(\tau^{(n)}_{g_i})_\ast ( [p_{g_i}])\neq 0, ~~~~~~(\tau^{(n)}_{g_{i_2}})_\ast ([p_{g_{i_1}}])=0~~~~when ~~i_1<i_2.$$
Now  Lemma 3.3 follows.
\qed

We are now ready to prove Theorem 1.4.

\noindent
Proof of Theorem 1.4:

By Theorem 2.2, it suffices to prove Conjecture 1.1 for groups coarsely embeddable into Hilbert space.

Let $\underline{E} G$ be the universal space for proper $G$-action.
We have the following commutative diagram:

$$\begin{array}[c]{ccc}
H_{0} ^{Or~G}(\underline{E} G, {\mathbb K}({\cal S})^{-\infty} ) &\stackrel{A}{\rightarrow}& K_{0}( {\cal S} G)\\
\downarrow\scriptstyle{\psi_\ast}&&\downarrow\scriptstyle{i_\ast}\\
K_{0} ^{G}(\underline{E}G) &\stackrel{\mu}{\rightarrow}& K_0( C^\ast(G))
\end{array},$$
where the right vertical map is induced by the inclusion map $i$ and the left vertical map is induced by the a map at the spectra level.

The strong Novikov conjecture holds  for groups coarsely embeddable into Hilbert space, i.e. the bottom horizontal map $\mu$ in the above diagram is injective \cite{Y1, STY}.
This, together with the fact that the left vertical map $\psi_\ast$ in the above diagram  is an isomorphism,  implies that the right vertical map $i_\ast$ in the above diagram is an isomorphism from the image of
$A: H_{0} ^{Or~G}(\underline{E} G, {\mathbb K}({\cal S})^{-\infty} ) \rightarrow   K_{0}( {\cal S} G)$,    to the image of $\mu:
  K_{0} ^{G}(\underline{E}G)  \rightarrow  K_0( C^\ast(G)) $.

We also have the  following commutative diagram:

$$\begin{array}[c]{ccc}
H_{0} ^{Or~G}(E G, {\mathbb K}({\cal S})^{-\infty} ) &\stackrel{A}{\rightarrow}& K_{0}( {\cal S} G)\\
\downarrow\scriptstyle{\psi_\ast}&&\downarrow\scriptstyle{i_\ast}\\
 K_{0} ^{G}(EG)  &\stackrel{\mu}{\rightarrow}& K_0( C^\ast(G))
\end{array},$$
where the left vertical map in the above diagram is the restriction of the left vertical map  in the previous diagram (here we identify  $H_{0} ^{Or~G}(E G, {\mathbb K}({\cal S})^{-\infty} )$ as a subgroup of
$H_{0} ^{Or~G}(\underline{E} G, {\mathbb K}({\cal S})^{-\infty} )$ and $K_{0} ^{G}(EG)$ as a subgroup of $ K_{0} ^{G}(\underline{E}G)   $).

 Assume that by contradiction  there is a nonzero element
$x$ in the abelian subgroup of $  K_0( C^\ast(G)) $ generated by
$\{[p_{g_1}], \cdots, [p_{g_n}]\}$ such that $x=\mu(y)$ for some $y\in K_{0} ^{G}(EG) .$
The left vertical map $\psi_\ast$ in the above diagram is  an isomorphism. This implies  $y= \psi_\ast (y')$ for some $y'\in H_{0} ^{Or~G}(E G, {\mathbb K}({\cal S})^{-\infty} )$. There exists an element  $x'$ in the finite part of $ K_{0}( {\cal S} G)$ such that
$i_\ast (x')=x$.
 We have $ i_\ast (A(y')-x')=0 .$  We observe that  the finite part of $ K_{0}( {\cal S} G)$ is  contained in the image of
 the map $A: H_{0} ^{Or~G}(\underline{E} G, {\mathbb K}({\cal S})^{-\infty} ) \rightarrow  K_{0}( {\cal S} G)$.
   As a consequence, we know that $A(y')-x'$ is in the image of the map $A: H_{0} ^{Or~G}(\underline{E} G, {\mathbb K}({\cal S})^{-\infty} ) \rightarrow  K_{0}( {\cal S} G)$. It follows from the isomorphism statement before the above diagram
  that $A(y')-x'=0$.
Hence $A(y')=x'$. Observe that $x'$ is a nonzero element in the finite part of $K_{0}( {\cal S} G)$.
This is a contradiction with Lemma 3.2. \qed

We remark that the argument in the above proof shows that the strong Novikov conjecture implies that
Conjecture 1.1.

\begin{theo}
Let $M$ be a compact oriented manifold with dimension $4k-1$ ($k>1$). If $\pi_1 (M )=G$ and $\{g_1, \cdots, g_n\}$ be a collection of finite order elements in $G$ such that $g_i\neq e$ for all $i$ and $\{[p_{g_1}], \cdots, [p_{g_n}]\}$ generates an abelian subgroup of $K_0(C^* (G))$ with rank $n$
and any nonzero element in the abelian subgroup of $K_0(C^* (G))$ generated by $\{[p_{g_1}], \cdots, [p_{g_n}]\}$ is not in the image of
of the map $\mu:
K_0^G(EG)\rightarrow K_0(C^*(G))$, then the rank of the structure group $S(M)$ is greater or equal to
$n$.
\end{theo}
\proof
We recall the surgery exact sequence:
$$\cdots \rightarrow H_{4k}(M, \mathbb{L} ) \rightarrow L_{4k}(\mathbb{Z} G) \rightarrow S(M) \rightarrow H_{4k-1}(M, \mathbb{L} )\rightarrow \cdots.$$
For each finite subgroup $H$ of $G$, we have the following commutative diagram:
$$\begin{array}[c]{ccc}
H_{4k} ^{H}(\underline{E}H , {\mathbb L}) &\stackrel{A}{\rightarrow}& L_{4k}( \mathbb{Z} H)\\
\downarrow&&\downarrow\\
H_{4k} ^{G}(\underline{E}G , {\mathbb L}) &\stackrel{A}{\rightarrow}& L_{4k}( \mathbb{Z} G)
\end{array},$$
where the vertical maps are induced by the inclusion homomorphism from $H$ to $G$.
For each element $g$ in $H$ with finite order, $p_g$ gives an element in $L_0 ( \mathbb{Q} H)$. Let
$[q_g]$ be the corresponding element in $L_{4k} ( \mathbb{Q} H)$ given by periodicity.
Recall that $$L_{4k} (\mathbb{Z} H) \otimes \mathbb{Q} \simeq L_{4k} (\mathbb{Q} H) \otimes \mathbb{Q}.$$
For each element $g$ in $H$ with finite order, we use the same notation $[q_g]$ to denote  the element in $L_{4k} (\mathbb{Z} H)\otimes \mathbb{Q}$ corresponding to   $[q_g]\in
L_{4k} ( \mathbb{Q} H)$ under the above isomorphism.

We have the following commutative diagram:
$$\begin{array}[c]{ccc}
H_{4k} ^{G}(EG , {\mathbb L})\otimes \mathbb{Q} &\stackrel{A}{\rightarrow}& L_{4k}( \mathbb{Z} G)\otimes \mathbb{Q} \\
\downarrow&&\downarrow\\
K_0^G (EG) \otimes \mathbb{Q} &   \stackrel{\mu}{\rightarrow}   & K_0( C^\ast (G))  \otimes \mathbb{Q}
\end{array},$$
where the left vertical map is induced by a map at the spectra level and the right vertical map is induced by the inclusion map:
$$ L_{4k}( \mathbb{Z} G)\rightarrow L_{4k}( C^\ast (G) ) \cong  K_0( C^\ast (G))$$ (see \cite{R2} for the last identification).

  The above commutative diagram, together with the assumption that any nonzero element in the abelian subgroup of $K_0(C^* (G))$ generated by $\{[p_{g_1}], \cdots, [p_{g_n}]\}$ is not in the image of
of the map $\mu:
K_0^G(EG)\rightarrow K_0(C^*(G))$, implies
 that (1) any nonzero element in the abelian subgroup of $L_{4k} ( \mathbb{Z} G)\otimes \mathbb{Q}$ generated by $\{[q_{g_1}], \cdots, [q_{g_n}]\}$
is not in the image of the rational assembly map $$A: H_{4k} ^{G}(EG , {\mathbb L}) \otimes \mathbb{Q} \rightarrow L_{4k} ( \mathbb{Z} G)\otimes \mathbb{Q};$$
(2) the abelian subgroup of $L_{4k} ( \mathbb{Z} G)\otimes \mathbb{Q}$ generated by $\{[q_{g_1}], \cdots, [q_{g_n}]\}$ has rank $n$.

By exactness of the surgery sequence,
we know that the the map: $L_{4k} ( \mathbb{Z} G)\otimes \mathbb{Q}\rightarrow S(M)\otimes \mathbb{Q},$ is injective on the abelian subgroup of  $L_{4k} ( \mathbb{Z} G)\otimes \mathbb{Q}$ generated by $\{[q_{g_1}], \cdots, [q_{g_n}]\}$.
This, together with the fact that the abelian subgroup of $L_{4k} ( \mathbb{Z} G)\otimes \mathbb{Q}$ generated by $\{[q_{g_1}], \cdots, [q_{g_n}]\}$ has rank $n$,  implies our theorem.
\qed

\begin{cor} Let $G$ be a countable group. If $G$ is finitely embeddable into Hilbert space and
 $M$ is a compact oriented manifold with dimension $4k-1$ $(k>1)$ and $\pi_1 (M)=G$,
then the rank of the structure group
$S(M)$ is greater than or equal to $N_{fin}(G)$.
\end{cor}

Let $G$ be a countable group. We define $r_{fin}(G)$ to be the rank of the abelian subgroup of  $K_0(C^* (G))$ generated by $[p_g]$
for all finite order elements $g$ in $G$. We emphasize that here we allow $g$ to be the identity $e$.

The proof of the following result is similar to that of Theorem 1.4 and is therefore omitted.
\begin{theo}
If $G$ is a countable group finiteley embeddable into Hilbert space, then  $r_{fin}(G)$  is greater than or equal to $N_{fin}(G)+1$.
\end{theo}

The following result is a consequence of Theorem 3.4 and part (1) of  Theorem 2.3.

\begin{cor} Let $M$ be a compact oriented manifold with dimension $4k-1$ $(k>1)$ and $\pi_1 (M)=G$.
 If there exists an element $g\neq e$ in $ G$ with finite order, then the rank of $S(M)$ is greater than or equal to one.
\end{cor}

The above corollary gives a different proof of a theorem of Chang and Weinberger \cite{CW}.

We mention that if the  group satisfies the strong Novikov conjecture, then the same method can be used to prove  a stronger statement than Theorem 3.4.
 In this case, one can show that  the number of conjugacy classes of  nontrivial  finite order elements in the group is a lower bound for both the rank of the finite part of
  operator K-theory and the rank of the structure group $S(M)$ when the dimension of $M$ is $4k-1$ $(k>1)$. To prove this result, we need to consider  K-theory classes induced by all representations of finite subgroups of the fundamental group.
However, our method doesn't yield the same lower bound for finitely embeddable groups.

In surgery theory, if $M$ is a compact oriented manifold,
the elements in the structure group $S(M)$ are pairs $(f, M')$, where $M'$ is a compact oriented manifold and
$f: M'\rightarrow M,$ is an orientation preserving  homotopy equivalence.
When an element is nontrivial, it is often the case that $M'$ and $M$ are homeomorphic, but $f$ is not  homotopy equivalent to a homeomorphism.
In \cite{CW}, Chang and Weinberger observed that the von Neumann trace can be used to distinguish the manifolds, not just structures.

We conjecture that elements of the structure group distinguished by the method of this paper are actually different manifolds.

Let  $M$ be a compact oriented manifold. Let $S_0(M)$ be the abelian subgroup of $S(M)$ generated by elements $[(f, M')]-[(\psi\circ f, M')]$, where $f: M'\rightarrow M$, is an orientation preserving homotopy equivalence and $\psi: M\rightarrow M,$ is  an orientation preserving  self homotopy equivalence.
We define the reduced structure group $ \tilde{S}(M)$ to be the quotient group $S(M)/S_0(M)$ (it is the coinvariant of the action of orientation preserving
self homotopy equivalence of $M$ on $S(M)$).

The following conjecture gives a lower bound on the ``size'' of the set of  different manifolds  in the structure group.

\begin{con}
If $M$ is a compact oriented manifold with dimension $4k-1$ $(k>1)$ and $\pi_1 (M)=G$,
then the rank of the reduced structure group
$\tilde{S}(M)$ is greater than or equal to $N_{fin}(G)$.
\end{con}

In a special case, we can verify this conjecture.

\begin{theo}  If $G$ has a homomorphism $\phi$ to a residually finite group such that $kernel (\phi)$ is torsion free,
then the above conjecture holds.
\end{theo}
\proof

For each integer $m\geq 1$, let $G_m$ be the intersection of all subgroups of $G$ with index at most $m$.
Observe that $G_m$ is a finite index subgroup of $G$ and $G_m$ is preserved under the action of $Aut(G)$, the group of all automorphisms of $G$.
As a consequence, the semidirect product $G\rtimes Aut(G)$ has a homomorphism  $\phi_m$ to the finite group $Q=G/G_m \rtimes Aut(G/G_m)$.

Let $\{g_1, \cdots, g_n\}$ be elements in $G$ with distinct finite orders $\{d_1, \cdots, d_n\}$. Without loss of generality, we can assume that $d_i< d_j$ when $i<j$.
Let $g_i'=\phi_m(g_i)$ in  $Q=G/G_m \rtimes Aut(G/G_m)$ for each $i$.
By the assumption of our theorem, there exists a sufficiently large integer $m$  such that $g_i'$  has order $d_i$ for all $1\leq i \leq n$. If $N$ is a compact oriented manifold with dimension $4k-1$ and $\pi_1(N)=G$,
there is an associated map: $\sigma: N \rightarrow BQ$ obtained as a composition of the classifying map from $N$ to $BG$ with the map induced from the inclusion map:
$G\rightarrow G\rtimes Aut(G),$ the homomorphism $\phi_m: G\rtimes Aut(G)\rightarrow Q$.
 The bordism group of such pairs $(N, \sigma)$ over $BQ$ is finite since $Q$ is finite and $4k-1$ is odd.
Hence there exists a positive integer $l$ such that $lN$ is the boundary of
 a compact oriented manifold  $W$ equipped with a map $h: W\rightarrow BQ$ that extends $\sigma$. Let $W_Q$ and $N_Q$ be respectively  the $Q$-covers  of $W$ and $N$.
 We note that $W_Q$ is the pull back of the $Q$-principal bundle of $EQ$ over $BQ$ by the map $h$.
  We define the following $g_i$-torsion invariant of $(N, \psi)$: $$\tau_{g_i}(N,\sigma)=\frac{1}{l} ( tr (g_i' |_{H_+})- tr(g_i' |_{H_{-}} )),$$
 where  $$H^{2k} (W_Q, N_Q) = H^{2k} _+(W_Q, N_Q)\oplus H^{2k}_{-}(W_Q, N_Q) $$
 is the decomposition corresponding to the positive and negative part of the symmetric bilinear form associated to the cup product on $H^{2k} (W_Q, N_Q)$,
 and $H_{+} = H^{2k} _+(W_Q, N_Q)$ and $H_{-}= H^{2k}_{-}(W_Q, N_Q) ,$
 and $tr$ is the standard trace on the algebra of  all linear operators on the finite dimensional vector space such that the trace of the identity is the dimension of the vector space. It is not difficult to see that this invariant is independent of the choice of $(W,\sigma)$.

 Note that if $ \psi$ is an orientation preserving  self homotopy equivalence of $N$, then we have  $$\tau_{g_i}(N,\sigma)=\tau_{g_i}(N,\sigma \circ \psi).$$
 The above identity can be shown as follows. Let $\psi_\ast$ be the element in $Aut(G)$ induced by the map $\psi$ and denote  the element   $\phi_m(\psi_\ast)$ in $Q$ by $\psi_\ast'$.
 Observe that the map $\sigma \circ \psi $ is homotopy equivalent to $\tilde{\psi}\circ \sigma $, where $\tilde{\psi}$
 is the map from $BQ$ to $BQ$ induced by the  homomorphism:  $Q\rightarrow Q$ defined by: $q\rightarrow \psi_\ast ' q (\psi_\ast ')^{-1}$.
 Let $\tilde{h}$ be the map from $W$ to $BQ$ defined by: $\tilde{h}=\tilde{\psi}\circ h$. It is easy to see that $\tilde{h}$ extends the map $\sigma \circ \psi$ from
 $N$ to $BQ$. By the definition of the $g_i$-torsion invariant of $(N, \sigma \circ \psi)$,
 we have  $$\tau_{g_i}(N,\sigma \circ \psi)=\frac{1}{l} ( tr ( \psi_\ast '  g_i'(\psi_\ast ')^{-1} |_{H_+})- tr(\psi_\ast ' g_i'  (\psi_\ast ')^{-1}|_{ H_{-}} )).$$
 Now our desired identity follows from the above equation and the trace property.

 The above identity implies  that the homomorphism from $S(M)$ to ${\mathbb R}$: $$(f, M') \rightarrow \tau_{g_i}(M', \theta \circ f  ),$$
 is $0$ on $S_{0}(M)$, where $\theta$ is the map from $M$ to $BQ$ obtained by  composing  the classifying map $M\rightarrow BG$ with the map from $BG$ to $BQ$  derived
 from the inclusion $G\rightarrow G\rtimes Aut(G)$ and the homomorphism $\phi_m$  from $G\rtimes Aut(G)$ to $Q$.
  As a consequence, this homomorphism from $S(M)$ to ${\mathbb R}$ induces a homomorphism from $\tilde{S}(M)$ to ${\mathbb R}$.

Let $q_{g_i}$ be the element in $L_{4k}({\mathbb Z} G)$ as in the proof of Theorem 3.4. There is a cobordism $Y_i$ realizing $q_{g_i}$ such that $\partial Y_i =M\cup (-M_i)$  for some compact manifold $M_i$ and the maps $$(f_i, \partial_0 f_i, \partial_1 f_i) : (Y_i, M, M_i) \rightarrow (M\times [0,1], M\times \{0\}, M\times \{1\})$$ satisfy the following conditions:

\noindent  (1) $\partial_0 f_i$, the restriction of $f_i$ to $M$, is the identity map from $M$ to $M$;

\noindent  (2) $\partial_1 f_i$, the restriction of $f_i$ to $M_i$,  is an orientation preserving homotopy equivalence from $M_i$ to $M$.

 We abbreviate $ \tau_{g_i}(M_j, \theta \circ (\partial_1 f_i) )$  by $ \tau_{g_i}(M_j )$ and $ \tau_{g_i}(M , \theta) $ by $  \tau_{g_i}(M).$
 We have
 $$ \tau_{g_i}(M_j )-\tau_{g_i}(M)=tr(g_i'p_{g_j'}),$$ where $$p_{g_j'}=\frac{1}{d_j}\sum_{k=1}^{d_j}(g_j')^k$$   and $tr$ is the canonical trace on the group algebra of $ Q$.
 It follows that
   $$\tau_{g_i}(M_j)-\tau_{g_i}(M)\neq 0~~~if ~~~i=j, ~~~\tau_{g_i}(M_j)-\tau_{g_i}(M)= 0~~~if ~~~i> j.$$
 This implies the abelian subgroup of $\tilde{S}(M)$ generated by $$[(f_1, M_1)]-[(id_M, M)], \cdots, [(f_n, M_n)] -[(id_M, M)]$$ has rank $n$,
 where $id_M$ is the identity map on $M$.
\qed

\begin{cor}
If $G$ is residually finite, then
Conjecture 3.8 holds.
\end{cor}

We remark that the proof of Theorem 3.9 shows that elements of the structure group distinguished by the method of this paper are indeed different manifolds in this special case.

\section{Applications to the  space of positive scalar curvature metrics}

In this section, we apply our result on the finite part of K-theory for group $C^\ast$-algebras  to estimate the size of the  space of positive scalar curvature metrics on a compact smooth spin manifold. This section is influenced by the previous work of Rosenberg and Stolz \cite{RS, S1, S2}.  In the finite group case, our result follows from the work of Stolz \cite{S1}.

If a compact smooth spin manifold $M$ has a positive scalar curvature metric and the dimension of $M$ is greater than or equal to 5, we introduce an abelian group  $P(M)$ of  equivalence classes of  all positive scalar curvature metrics on $M$.
 We  give a lower bound on the rank of $P(M)$ when the dimension of $M$ is $2k-1$ with $k>2$ and the fundamental group $\pi_1(M)$ is finitely embeddable into Hilbert space. For general group, we show that if $\pi_1 (M)$ is not torsion free, then the rank of
the abelian group $P(M)$ is at least one when dimension of $M$ is $2k-1$ $(k>2)$ and is at least two when dimension of $M$ is $4k-1$ $(k>1)$ .
We remark that we don't need the spin condition to define $P(M)$, but we  require the spin condition to get a lower bound on its size.

Let $M$ be a compact smooth manifold $M$  with $\pi_1 (M)=G$ and dimension greater than or equal to five. Assume that $M$ has a positive scalar curvature metric $g_M$.
Let $I$ be the closed interval.
We first form a connected sum $(M\times I) \sharp (M\times I)$, where the connected sum is performed away from the boundary of each copy of $M\times I$. We define the generalized connected sum $(M\times I) \natural (M\times I)$ to be  the manifold obtained from $(M\times I) \sharp (M\times I)$  by surgering the kernel of the homomorphism from  $\pi_1 ((M\times I) \sharp (M\times I)) =G\ast G$ to $G$. $(M\times I) \natural (M\times I)$  has four boundary components, two components being $M$ and the other two being $-M$, where $-M$ is the manifold $M$ with reversed orientation.
 If $g_1$ and $g_2$ be positive scalar curvature metrics on $M$, we place $g_M$ on one boundary component $M$ and $g_1$ and $g_2$ on the other two boundary components $-M$. By the surgery theorem \cite{GL, SY},
 there exists a positive scalar curvature metric on $(M\times I) \natural (M\times I)$
 such that it is a product metric near the other boundary component $M$. We denote by  $g$ the positive scalar curvature metric on this boundary component $M$.

 We prove that if $g$ and $g'$ are two positive scalar curvature metrics $g$ and $g'$ on $M$  obtained from the same pair of positive scalar curvature metrics $g_1$ and $g_2$ by the above process. We prove that $g$ and $g'$ are concordant. We can glue the remaining boundaries of the two generalized connected sums to form a cobordism of $M$
 to  itself (the metric  is a product near the boundary and has positive scalar curvature) with the fundamental group surjecting into $G$.
 We surger away the kernel of the surjection.  We note that the nullcobordisms of a given $n$-manifold $N$ with a homomorphism of $\pi_1(N)\rightarrow G$ is a torsor,
 i.e. given one nullcobordism, the remaining nullcobordisms, up to nullcorbordism relative to the boundary, form an abelian group $( \cong \Omega_{n+1}( N) )$.
 Here the preferred null cobordism is $N \times I$. We claim that the element given by our cobordism is trivial--this is because  we can glue the boundaries of the cylinders $M\times I$ in the cobordism  to obtain $(M\times S^1) \natural (M\times S^1)$ which bounds $(M\times D^2) \natural (M\times D^2) $.
 Consequently we can apply the surgery theorem   \cite{GL, SY} to this cobordism of $(M, g)$ to $(M, g')$ to obtain a new positive scalar curvature metric cobordism to these manifolds
 where the underlying manifold is $M\times I$. Thus $g$ and $g'$ are concordant.

 Two positive scalar curvature metrics $g$ and $g'$ on $M$ are defined to be
 equivalent if and only if they are concordant.
 We define the equivalence class $[g]$ to be the sum of the equivalence classes of $[g_1]$ and $[g_2]$ with respect to $[g_M]$.
 By argument similar to that in the previous paragraph, we can show that the set of the concordance classes of all positive scalar curvature metrics on $M$ is an abelian semigroup with respect to this sum operation.
 We define the abelian group $P(M)$ to be the Grothendick group of the above abelian semigroup. We point put that our definition of $P(M)$ is closely related to Stolz's group $R_n$ of concordance classes of positive scalar curvature metrics \cite{S2}.

 Recall that $r_{fin}(G)$ is the rank of the abelian subgroup  of $K_0( C^*(G))$ generated by $[p_g]$ for all finite order elements $g$ in $G$.
 Here $g$ is allowed to be the identity element $e$.

\begin{theo} \noindent (1) Let $M$ be  a compact smooth spin manifold with a positive scalar curvature metric and dimension  $2k-1$ ($k>2$).
The rank of the abelian group $P(M)$ is greater than or equal to $r_{fin}(G)-1$;

\noindent (2) Let $M$ be  a compact smooth spin manifold with a positive scalar curvature metric and dimension  $4k-1$ ($k>1$).
The rank of the abelian group $P(M)$ is greater than or equal to $r_{fin}(G)$.
\end{theo}

The following result is a consequence of the above theorem and Theorem 3.6.

\begin{cor}\noindent (1) Let $M$ be  a compact smooth spin manifold with a positive scalar curvature metric and  dimension  $2k-1$ ($k>2$).
If $\pi_1(M)=G$ is finitely embeddable into Hilbert space, then the rank of the abelian group $P(M)$ is greater than or equal to $N_{fin}(G);$

\noindent (2) Let $M$ be  a compact smooth spin manifold with a positive scalar curvature metric and the dimension  $4k-1$ ($k>1$).
If $\pi_1(M)=G$ is finitely embeddable into Hilbert space, then the rank of the abelian group $P(M)$ is greater than or equal to $N_{fin}(G)+1.$
\end{cor}

The following result is a  consequence of Theorem 4.1 and part (2) of Theorem 2.3.

\begin{cor}
\noindent (1) Let $M$ be  a compact smooth spin manifold with a positive scalar curvature metric and  dimension  $2k-1$ ($k>2$).
 If $\pi_1(M)$ is not torsion free, then the rank of the abelian group $P(M)$ is greater than or equal to one;

\noindent (2) Let $M$ be  a compact smooth spin manifold with a positive scalar curvature metric and the dimension  $4k-1$ ($k>1$).
 If $\pi_1(M)$ is not torsion free, then the rank of the abelian group $P(M)$ is greater than or equal to two.
\end{cor}

 Piazza and Schick used a different method to prove that the  space of positive scalar curvature has infinitely many connected components when
$M$ is  a compact smooth spin manifold with a positive scalar curvature metric and dimension  $2k-1$ ($k>2$), and
the fundamental group $\pi_1(M)$ is not torsion free \cite{PS}.

We need some preparation to prove the main result in this section.
Let $F$ be a finite group and $N$ be a $F$-manifold. We say that $N$ is $F$-connected if $N/F$ is connected.

\begin{prop}
For given positive integers $d$ and  $k$,

\noindent (1) there exist ${\mathbb Z}/d{\mathbb Z}$-connected  compact smooth spin ${\mathbb Z}/d{\mathbb Z}$-manifolds $\{N_1, \cdots, N_n\}$
 such that
the dimension of each $N_i$ is $2k$ and
the ${\mathbb Z}/d{\mathbb Z}$-equivariant indices of the Dirac operators on $\{N_1, \cdots, N_n\}$ rationally generate  $RO({\mathbb Z}/d{\mathbb Z})\otimes {\mathbb Q}$ modulo rational multiples of  the regular representation, the ${\mathbb Z}/d{\mathbb Z}$ action on each $N_l$ ($1\leq l\leq n$) is free except for finitely many fixed points of
the ${\mathbb Z}/d{\mathbb Z}$ action, where   $RO({\mathbb Z}/d{\mathbb Z})$ is  the real representation ring of $  {\mathbb Z}/d{\mathbb Z}$;

\noindent  (2) there exist ${\mathbb Z}/d{\mathbb Z}$-connected  compact smooth spin ${\mathbb Z}/d{\mathbb Z}$-manifolds $\{N_1, \cdots, N_n\}$
 such that
the dimension of each $N_i$ is $4k$ and
the ${\mathbb Z}/d{\mathbb Z}$-equivariant indices of the Dirac operators on $\{N_1, \cdots, N_n\}$  rationally generate  $RO({\mathbb Z}/d{\mathbb Z})\otimes {\mathbb Q}$,
 and the ${\mathbb Z}/d{\mathbb Z}$ action on each $N_l$ ($1\leq l\leq n$) is free except for finitely many fixed points of
the ${\mathbb Z}/d{\mathbb Z}$ action.
\end{prop}
\proof
(1)
Let $\{a_1, \cdots, a_d\}$ be distinct points in ${\mathbb C}-\{0\}.$
Define a two dimensional smooth compact surface $S$ by:
$$ S=\{ [x,y,z]\in {\mathbb C}{\mathbb P}^2: y^d=(x-a_1 z) \cdots (x-a_d z)\}.$$
Observe that $S$ is a surface with genus $\frac{(d-1)(d-2)}{2}$.
We have

\noindent (a) $[0,1,0]$ is not in $S$;

\noindent (b) the map $\pi: S\rightarrow {\mathbb C}{\mathbb P} ^1$ sending $[x,y,z]$ to $[x,z]$,  is a  branched covering with degree $d$;

\noindent (c) there exist $d$ number of branch points $\{ [a_1, 1], \cdots, [a_d, 1]  \} \subset {\mathbb C}{\mathbb P}^1.$

For each positive integer $1\leq l \leq d$, we have a natural action ${\mathbb Z}/d{\mathbb Z}$ on $S$ by:
$$ \alpha_l([m]) [x,y,z]= [x, exp(\frac{2\pi i ml}{d})y, z]$$ for any $[m]\in {\mathbb Z}/d{\mathbb Z}$ and $[x,y,z]\in S$.

We shall prove Part (1) of Proposition 4.4 by induction.
We first deal with the special case that  $d$ is $2$ or a prime number.

Let $N_l$ be  the product of  $k$ number of copies  of $S$ with  the diagonal action of $ \alpha_l$. The assumption that $d$ is $2$ or a prime number
implies that  the action of ${\mathbb Z}/d{\mathbb Z}$
on each $N_{l}$ is free  except for finitely many fixed points.

By Theorem 8.35 in Atiyah-Bott \cite{AB}, we know that the  ${\mathbb Z}/d{\mathbb Z}$-equivariant indices of the Dirac operators on $\{N_1, \cdots, N_d\}$ generate $RO(  {\mathbb Z}/d{\mathbb Z})\otimes {\mathbb Q}$ modulo rational multiples of  the regular representation.

Assume by induction that the proposition is true for any divisor of $d$ less than $d$.
Let $\{N_1, \cdots N_s\}$ be the set of ${\mathbb Z}/d{\mathbb Z}$-connected  compact smooth spin ${\mathbb Z}/d{\mathbb Z}$-manifolds
whose ${\mathbb Z}/d{\mathbb Z}$-equivariant indices of Dirac operators rationally generate all real  representations induced from all proper subgroups modulo rational multiples of the regular representation.
It suffices to construct a set of ${\mathbb Z}/d{\mathbb Z}$-connected  compact smooth spin ${\mathbb Z}/d{\mathbb Z}$-manifolds with dimension $2k$
whose ${\mathbb Z}/d{\mathbb Z}$-equivariant indices of Dirac operators rationally generate all real representations with
 characters supported on  each generator of the abelian group ${\mathbb Z}/d{\mathbb Z}$ and its inverse.

Let $\{[ l_1], \cdots, [l_m]\}$ be the set of all generators of the abelian group ${\mathbb Z}/d{\mathbb Z}$.
Let $N_{s+l_i}$ be  the product of  $k$ number of copies  of $S$ with  the diagonal action of $ \alpha_{l_i}$.
The fact that $l_i$ is a generator of the abelian group ${\mathbb Z}/d{\mathbb Z}$ implies that the action of ${\mathbb Z}/d{\mathbb Z}$
on each $N_{s+l_i}$ is free  except for finitely many fixed points.
By the induction hypothesis and Theorem 8.35 in Atiyah-Bott \cite{AB}, we know that the  ${\mathbb Z}/d{\mathbb Z}$-equivariant indices of the Dirac operators on $\{ N_{s+l_1}, \cdots, N_{s+l_m}\}$ rationally generate  all real representations with
 characters supported on each generator of the abelian group ${\mathbb Z}/d{\mathbb Z}$ and its inverse modulo the abelian subgroup rationally generated by all
 representations induced from proper subgroups of ${\mathbb Z}/d{\mathbb Z}$. This, together with the induction hypothesis, implies  that ${\mathbb Z}/d{\mathbb Z}$-equivariant indices of the Dirac operators on
 $\{N_1, \cdots, N_s,    N_{s+l_1}, \cdots, N_{s+l_m}\}$   generate $RO(  {\mathbb Z}/d{\mathbb Z})\otimes {\mathbb Q}$ modulo rational multiples of the regular representation.

 (2) For the second part of the proposition, when $d=1$, ${\mathbb Z}/d{\mathbb Z}$ is a trivial group and we can take a $4k$-dimensional compact smooth spin manifold whose Dirac operator has nonzero index. For example, we can take the product of $k$ copies of the Kummer (or K3) surface:
 $$\{ (z_0, z_1, z_2, z_3): z_0^4 +z_1^4 +z_2^3+z_3^4=0\} \subseteq {\mathbb  C} {\mathbb P}^3.$$
 This proves part (2) for the trivial group case.
 The rest of the proof goes exactly the same as in the proof of of part (1) with the dimension of the manifold changed to $4k$. We point out that in the induction process, the trivial subgroup  induces the regular representation.
\qed

We remark that Proposition 4.4 can be generalized to any finite group by using representations induced from its cyclic subgroups.

Now we are ready to prove Theorem 4.1.
\noindent \proof
Let $\tilde{M}$ be the universal cover of $M$. We will prove Theorem 4.1 when the dimension of $M$ is $4k-1$. The proof is completely similar when the dimension is $2k-1$.

For each finite order element $g$ in $G$ with order $d$. By Proposition 4.4, there exist ${\mathbb Z}/d{\mathbb Z}$-connected  compact smooth spin ${\mathbb Z}/d{\mathbb Z}$-manifolds $\{N_1, \cdots, N_n\}$
 such that
the dimension of each $N_i$ is $4k$ and the sum of
the ${\mathbb Z}/d{\mathbb Z}$-equivariant indices of the Dirac operators on  $\{ N_1, \cdots, N_n\}$ is a nonzero multiple of
the trivial representation of ${\mathbb Z}/d{\mathbb Z}$.

Let $N_{g, l} = G\times_{{\mathbb Z}/d{\mathbb Z}}N_l, $ where ${\mathbb Z}/d{\mathbb Z} $ acts on $N_l$ as in Proposition 4.4 and ${\mathbb Z}/d{\mathbb Z}$ acts on $G$ by $[m]h=hg^m$ for all $h\in G$ and $[m]\in {\mathbb Z}/d{\mathbb Z}  $. Observe that $N_{g,l}$ is a $G$-manifold.

Let $\{g_1, \cdots, g_r\}$ be a collection of finite order elements such that \linebreak
 $\{[p_{g_1}], \cdots, [p_{g_r}]\}$ generates an abelian subgroup of
$K_0(C^\ast(G))$ with rank $r$. Let $N_{g_i}=\bigsqcup_{l=1}^{j_i}N_{g_i, l}$ be the disjoint union of all $G$-manifolds described as above. Let $I$ be the unit interval
$[0,1]$.
We first form a generalized  $G$-equivariant connected sum $(\tilde{M}\times I)\natural N_{g_i}$ along a free $G$-orbit of each  $N_{g_i,l}$ and away from the boundary of $\tilde{M}\times I$ as follows.
We first obtain a $G^{\ast j_i}$-equivariant  connected sum $(\tilde{M}\times I) \sharp N_{g_i}$ along a free $G$-orbit of each  $N_{g_i, l}$ and away from the boundary of $\tilde{M}\times I$, where $G^{\ast j_i}$ is the free product of $j_i$ copies of $G$.  More precisely,  we inductively form the $G^{\ast j_i}$-equivariant connected sum
$( \cdots ((\tilde{M}\times I) \sharp N_{g_i,1})\cdots )\sharp N_{g_i, j_i}$, where the equivariant connected sum is inductively taken along a free orbit
and away from the boundary.
We denote this space by $(\tilde{M}\times I) \sharp N_{g_i}$.
  We then perform surgeries on  $(\tilde{M}\times I) \sharp N_{g_i}$
 to obtain a $G$-equivariant cobordism between two copies of $G$-manifold $\tilde{M}$.

For any positive scalar curvature metric $h$ on $M$, by Theorem 2.2 \cite{RW},
the above cobordism gives us another positive scalar curvature metric $h_i$ on $M$.
By the assumption on $\{[p_{g_1}], \cdots, [p_{g_r}]\}$, the relative higher index theorem \cite{XY}, the relative higher index of the Dirac operator $M\times {\mathbb R}$ associated to the positive scalar curvature  metrics of  $h_i$ and $g_M$ is $[p_{g_i}]$ in $K_0(C^\ast(G))$. As a consequence,   we know that
$\{ [h_1], \cdots, [h_r]\}$ generates an abelian subgroup of $P(M)$ with  rank $r$.
\qed

We mention that if the fundamental group satisfies the strong Novikov conjecture, then the same method can be used to prove that the rank of $P(M)$
is greater than or equal to the number of conjugacy classes of nontrivial finite order elements in the fundamental group when the dimension of $M$ is $2k-1$ $(k>2)$,
and the rank of $P(M)$
is greater than or equal to the number of conjugacy classes of finite order elements in the fundamental group when the dimension of $M$ is $4k-1$ $(k>1)$.
To prove this result, we need to consider  K-theory classes induced by all representations of finite subgroups of the fundamental group.

\section{Groups finitely embeddable into Hilbert space}

In this section, we show that various classes of groups are finitely embeddable into Hilbert space.
In particular, we show that this class of groups include all residually finite groups, amenable groups, Gromov's monster groups, all
virtually torsion free groups (e.g. $Out(F_n)$), and all groups of analytic diffeomorphisms of an analytic connected  manifold fixing a given point.

We first introduce a concept of groups  locally embeddable into Hilbert space.

\begin{defi}  A countable discrete group $G$ is said to be locally embeddable into Hilbert space $H$ if for any finite subset $F\subseteq G$,
there exists a group $G'$ coarsely embeddable into $H$ such that there is a map $\phi: F\rightarrow G'$ satisfying
(1) $\phi(e)=e$ if $e\in F$;
(2) $\phi(gh)=\phi(g) \phi (h)$ if $g, h\in F$ and $gh\in F$;
(3) if $g$ and $h$ are distinct elements in $F$, then $\phi(g)\neq \phi(h)$.

\end{defi}

Clearly groups locally embeddable into Hilbert spaces are finitely embeddable into Hilbert space.
Recall that a group $G$ is said to be locally embeddable into finite groups (LEF) if the group $G'$ in the above definition can always to be chosen to be a finite group. Observe that all residually finite groups are LEF.

\begin{prop}
Let $N$ be an analytic connected  manifold and $x_0\in N$. If $G$ is a countable group of analytic diffeomorphisms fixing the point $x_0$,
then $G$ is locally embeddable into Hilbert space.
\end{prop}
\proof For any positive integer $k$, let $J_k$ be the finite dimensional vector space of all $k$-th jets at $x_0$.
Let $GL(J_k)$ be the Lie group of all linear isomorphisms  from $J_k$ to $J_k$.
Any diffeomorphism of $N$ fixing $x_0$ induces an isomorphism of $J_k$.  It follows that, for each $k$, there is a natural homomorphism $\psi_k$
from $G$ to $GL(J_k)$. Let $G'= \psi_k (G)$. $G'$ is coarsely embeddable into Hilbert space since $GL(J_k)$ is a Lie group with finitely many connected components \cite{GHW}. For any finite subset $F$ of $G$, we restrict $\psi_k$ to $F$ to obtain a map $\phi$ from $F$ to $G'$.
By analyticity, we can verify that $G'$ and $\phi$ satisfy the conditions in Definition 5.1 when $k$ is large enough.
\qed

It is an open question whether an arbitrary countable group of analytic diffeomorphisms on an analytic connected manifold is locally or finitely embeddable into Hilbert space. More generally, it is an open question whether an arbitrary countable group of diffeomorphisms on a smooth connected manifold is locally or finitely embeddable into Hilbert space.

The next result states that certain limit groups are locally embeddable into Hilbert space.

\begin{prop}
If $G$ is the limit of a sequence of groups $\{G_k\}_{k=1}^\infty$ coarsely embeddable into Hilbert space such that, for each finite subset $F\subset G_k$ for some $k$, there exists $l\geq k$ such that the map from the image of $F$ in $G_l$ to $G$ is injective,  then $G$ is locally embeddable into Hilbert space.
\end{prop}

We point out that in the above proposition  the homomorphism from $G_k$ to $G_l$ ($k>l$) is not assumed to be injective.
Examples of groups satisfying the above proposition include Burnside groups and Gromov's monster groups \cite{G, AD}.
The proof of the above result is straightforward and is therefore omitted.

\begin{prop}
 If $G$  has a   torsion free normal subgroup $G_0$ such that $G/G_0$ is residually finite,
then $G$ is finitely embeddable into Hilbert space.
\end{prop}
\proof
Let $\phi$ be the quotient homomorphism from $G$ to $G'=G/G_0$. We have $order (\phi(g))=order(g)$ for all finite order elements
$g$ in $G$.  This can be seen as follows. The order of $\phi(g)$ is a divisor of the order of $g$. If by contradiction that
$g$ is a finite order element in $G$ whose order is greater than $1$ and
$order (\phi(g))\neq order(g)$, then $order(g)=n ~~order (\phi(g))$ for some positive integer $n>1$. This would imply that $g^n$ is in the kernel of $\phi$. Hence $g^n$ is in $G_0$. But $g^n$ is a nontrivial finite order element. This is a contradiction with the assumption that $G_0$ is torsion free.
Now our proposition follows from the definition of residually finiteness of $G'$. \qed

\begin{cor}
If $G$ is virtually torsion free, then $G$ is finitely embeddable into Hilbert space.
\end{cor}

An example of a group satisfying the assumption in the above corollary is $Out(F_n)$ (Proposition 1.2 in \cite{B}).

Sofic groups have many common examples with finitely embeddable groups. It is an open question if sofic groups are finitely embeddable into Hilbert space.
Finally, we mention that it is also an open question to construct a countable group not finitely embeddable into Hilbert space.

\section{Appendix}

In this appendix, we prove the following analogue of Atiyah's $L^2$-index theorem for the maximal group $C^\ast$-algebra.
This result is a folklore. For completeness, we include a proof here.

\begin{theo} Let $G$ be a countable group.
Let $X$ be a a locally compact space with a proper, free, and cocompact action of $G$.
Let $\mu$ be the assembly map: $ K_0^G(X)  \rightarrow K_0(C^\ast (G))$.
If $\tau$ is the trace on $C^\ast (G)$ given by $\tau(\sum c_g g)=\sum c_g$, then $\tau (\mu [(F,H)])= index (F)$
for any $K$-homology class $[(F, H)]$ in  $K_0 ^G(X),$ where $index (F)$ is the Fredholm index of $F$.
\end{theo}
\proof
In  \cite{WY}, we give a different proof of Atiyah's $L^2$-index theorem. We shall use a key ingredient in this proof.

Let $d$ be a $G$-invariant metric on $X$ (compatible with the topology of $X$).
The assumption that the $G$ action on $X$ is proper and free implies that  there exists $\delta>0$ such that $d(x, gx)\geq 10\delta$ for all $x\in X$ and $g\neq G$ in $G$.
By the proof of the Atiyah's $L^2$-index theorem  in \cite{WY} (page 1402 of \cite{WY}),
we can write
$$\mu [(F, H)]=[p]- \bigl(\begin{smallmatrix}
1&0\\ 0&0
\end{smallmatrix} \bigr),
$$
where
$$p
=p_1+ \bigl(\begin{smallmatrix}
1&0\\ 0&0
\end{smallmatrix} \bigr) $$
such that $p_1$ is an element in the matrix algebra over $({\cal S}_1 X)^G$,
the algebra of $G$-invariant and locally traceable operators
on a $G$-$X$-module $H$ with finite propagation,  satisfying $$propagation(p_1)<\delta.$$

Note that $({\cal S}_1 X)^G$ is isomorphic to ${\cal S}_1 G$, the group algebra of $G$ over the ring ${\cal S}_1$.
We can extend the trace $\tau$ naturally to a trace on
We have $$\tau (\mu [(F,H)])=\tau (p_1)=\sum_{g\in G} tr(g^{-1} p_1),$$
where $tr$ is the canonical trace on ${\cal S}_1 G$ as in Lemma 3.1.
We also note that the sum in the above identity is a finite sum since $p_1$ has finite propagation.
Now by Lemma 3.1 and the choice of $\delta$, we have $tr(g^{-1} p_1)=0$ for all $g\neq e$.
As a consequence, we have $\tau(p_1)=tr(p_1).$ By Atiyah's $L^2$-index theorem \cite{A}, we have $tr(p_1)= index (F)$.
It follows that $\tau (\mu [(F,H)])= index (F)$.
\qed

\noindent SW: Department of Mathematics, University of Chicago,

\noindent 5734 S. University Avenue
Chicago, IL 60637, USA.

\noindent e-mail: shmuel@math.uchicago.edu

\noindent GY: Department of Mathematics, Texas A\&M University,

\noindent College Station, TX 77843, USA, and

\noindent
Shanghai Center for Mathematical Sciences, China.

\noindent e-mail: guoliangyu@math.tamu.edu

\end{document}